\documentclass[twoside]{amsart}

\usepackage{latexsym}
\usepackage{amsthm,amsfonts,amssymb}
\usepackage[leqno]{amsmath}
\usepackage[all]{xy} \SelectTips{eu}{}
\usepackage{hyperref}

\newcommand{\numberseries}{\bfseries}   

\newlength{\thmtopspace}                
\newlength{\thmbotspace}                
\newlength{\thmheadspace}               
\newlength{\thmindent}                  

\newtheoremstyle{bfupright head,slanted body}
{\thmtopspace}{\thmbotspace}
{\itshape}{\thmindent}{\bfseries}{.}{\thmheadspace}
{{\numberseries \thmnumber{#2.\;}}\thmnote{#3}}

\newtheoremstyle{bfupright head,upright body}
{\thmtopspace}{\thmbotspace}
{\upshape}{\thmindent}{\bfseries}{.}{\thmheadspace}
{{\numberseries \thmnumber{#2.\;}}\thmnote{#3}}

\newtheoremstyle{it head,upright body}
{\thmtopspace}{\thmbotspace}
{\upshape}{\thmindent}{\upshape}{}{\thmheadspace}
{{\numberseries\thmnumber{#2.\!\negthickspace}}
  {}}

\newtheoremstyle{fixed bf head,slanted body}
{\thmtopspace}{\thmbotspace}{\itshape}
{\thmindent}{\bfseries}{.}{\thmheadspace}
{{\numberseries \thmnumber{#2.\;}}\thmname{#1}\thmnote{ (#3)}}

\newtheoremstyle{fixed bf head,upright body}
{\thmtopspace}{\thmbotspace}{\upshape}
{\thmindent}{\bfseries}{.}{\thmheadspace}
{{\numberseries \thmnumber{#2.\;}}\thmname{#1}\thmnote{ (#3)}}

\newtheoremstyle{independent paragraph}
{\thmtopspace}{\thmbotspace}
{\upshape}{\parindent}{\upshape}{}{0pt}
{\thmnote{#3 }}

\newtheoremstyle{subparagraph}
{\thmbotspace}{\thmbotspace}
{\upshape}{\parindent}{\upshape}{}{0pt}
{\thmnote{#3 }}

\theoremstyle{bfupright head,slanted body}
\newtheorem{res}{}[section]             \newtheorem*{res*}{}

\theoremstyle{bfupright head,upright body}
\newtheorem{bfhpg}[res]{}               \newtheorem*{bfhpg*}{}

\theoremstyle{it head,upright body}
                 \newtheorem*{com*}{}

\theoremstyle{fixed bf head,slanted body}
\newtheorem{thm}[res]{Theorem}          \newtheorem*{thm*}{Theorem}
      \newtheorem*{prp*}{Proposition}
\newtheorem{cor}[res]{Corollary}        \newtheorem*{cor*}{Corollary}
\newtheorem{lem}[res]{Lemma}            \newtheorem*{lem*}{Lemma}

\theoremstyle{fixed bf head,upright body}
      \newtheorem*{obs*}{Observation}
\newtheorem{rmk}[res]{Remark}           \newtheorem*{rmk*}{Remark}

\theoremstyle{independent paragraph}
\newtheorem{ipg}{}

\theoremstyle{subparagraph}

\newlength{\thmlistleft}        
\newlength{\thmlistright}       
\newlength{\thmlistpartopsep}   
\newlength{\thmlisttopsep}      
\newlength{\thmlistparsep}      
\newlength{\thmlistitemsep}     

\newcounter{eqc} 
\newenvironment{eqc}{\begin{list}{\upshape (\textit{\roman{eqc}})}%
    {\usecounter{eqc}%
      \setlength{\leftmargin}{\thmlistleft}%
      \setlength{\labelwidth}{\thmlistleft}%
      \setlength{\rightmargin}{\thmlistright}%
      \setlength{\partopsep}{\thmlistpartopsep}%
      \setlength{\topsep}{\thmlisttopsep}%
      \setlength{\parsep}{\thmlistparsep}%
      \setlength{\itemsep}{\thmlistitemsep}}}%
  {\end{list}}%

\newcommand{\eqclbl}[1]{{\upshape(\textit{#1})}}

\newcounter{prt}
  {\end{list}}%

\newcounter{rqm}
  {\end{list}}%

  {\end{list}}%

  {\end{list}}%

\newenvironment{prf}[1][Proof]{\begin{proof}[\bf #1]}{\end{proof}}

\newcommand{\pgref}[1]{(\ref{#1})}
\renewcommand{\eqref}[1]{\pgref{eq:#1}}

\renewcommand{\theequation}{\arabic{equation}}
\numberwithin{equation}{res}

\newcommand{\Hom}[3][R]{\operatorname{Hom}_{#1}(#2,#3)}
\newcommand{\tp}[3][R]{#2\otimes_{#1}#3}

\newcommand{\tpP}[3][R]{(\tp[#1]{#2}{#3})}

\newcommand{\m}{\mathfrak{m}}

\renewcommand{\H}[2][\no]{\operatorname{H}_{#1}(#2)}
\newcommand{\xra}{\xrightarrow}
\newcommand{\xla}{\xleftarrow}

\newcommand{\xri}{\xra{\;\is\;}}
\newcommand{\xle}{\xla{\;\eq\;}}
\newcommand{\xre}{\xra{\;\eq\;}}
\newcommand{\eq}{\simeq}
\newcommand{\is}{\cong}

\renewcommand{\le}{\leqslant}

\newcommand{\Coker}[1]{\mbox{\ensuremath{\operatorname{Coker}#1}}}
\newcommand{\mapdef}[4][\rightarrow]{\mbox{\ensuremath{#2\colon #3 #1 #4}}}
\newcommand{\dmapdef}[4][\lora]{#2\colon #3\:#1\:#4}

\newcommand{\Cone}{\operatorname{Cone}}

\newcommand{\ZZ}{\mathbb{Z}}

\newcommand{\Rmk}{(R,\m,k)}

\newcommand{\E}[2][R]{\operatorname{E}_{#1}(#2)}
\newcommand{\Kac}{\operatorname{{\bf K}_{ac}(Inj\ R)}}
\newcommand{\Ktac}{\operatorname{{\bf K}_{tac}(Inj\ R)}}
\newcommand{\Kinj}{\operatorname{{\bf K}(Inj\ R)}}
\newcommand{\Kcinj}{\operatorname{{\bf K}^c(Inj\ R)}}
\newcommand{\Df}{\operatorname{{\bf D}}^f(R)}
\newcommand{\Loc}{\operatorname{Loc}}

\newcommand{\lora}{\longrightarrow}
\newcommand{\Susp}[2][]{{\sf\Sigma}^{#1}{#2}}
\newcommand{\iR}{\mathsf{i}R}

\begin{document}

\allowdisplaybreaks[4]

\title{A test complex for Gorensteinness}



\author[L.~W.~Christensen]{Lars Winther Christensen}

\thanks{L.W.C.\ was partly supported by a grant from The Carlsberg
  Foundation.}

\address{Lars Winther Christensen, Department of Mathematics,
  University of Nebraska, Lincoln, Nebraska~68588}

\email{winther@math.unl.edu}

\urladdr{http://www.math.unl.edu/{\tiny $\sim$}lchristensen}

\author[O.~Veliche]{Oana~Veliche}

\address{Oana Veliche, Department of Mathematics, University of Utah,
  Salt Lake City, Utah~84112}

\email{oveliche@math.utah.edu}

\urladdr{http://www.math.utah.edu/{\tiny $\sim$}oveliche}

\date{17 January 2007}

\begin{abstract}
  Let $R$ be a commutative noetherian ring with a dualizing complex.
  By recent work of Iyengar and Krause \cite{SInHKr06}, the difference
  between the category of acyclic complexes and its subcategory of
  totally acyclic complexes measures how far $R$ is from being
  Gorenstein. In particular, $R$ is Gorenstein if and only if every
  acyclic complex is totally acyclic.
  
  In this note we exhibit a specific acyclic complex with the property
  that it is totally acyclic if and only if $R$ is Gorenstein.
\end{abstract}

\maketitle

\thispagestyle{empty}
\vspace*{-2ex}
\enlargethispage{3ex}
\section*{Introduction}

Let $R$ be a commutative noetherian ring. A complex $X$ of $R$-modules
is said to be \emph{acyclic} if it has zero homology, i.e.\ $\H{X}
=0$. An acyclic complex of projective modules is called \emph{totally
  acyclic} if the acyclicity is preserved by $\Hom{-}{P}$ for every
projective module $P$.  Dually, an acyclic complex of injective
modules is totally acyclic if the acyclicity is preserved by
$\Hom{I}{-}$ for every injective module~$I$.

Over a Gorenstein ring, every acyclic complex of projective or of
injective modules is totally acyclic. Iyengar and Krause have recently
proved a converse; indeed, by \cite[cor.~5.5]{SInHKr06} the following
are equivalent when $R$ has a dualizing complex:
\begin{eqc}
\item The ring $R$ is Gorenstein.
\item Every acyclic complex of projective $R$-modules is totally
  acyclic.
\item Every acyclic complex of injective $R$-modules is totally
  acyclic.
\end{eqc}
Moreover, for a local ring $(R,\m)$ that is not Gorenstein and has
$\m^2=0$ there is a natural example, provided by
\cite[prop.~6.1(3)]{SInHKr06}, of an acyclic, but not totally acyclic,
complex of projective $R$-modules.

\begin{ipg}
  The purpose of this note is to prove that for every ring $R$ with a
  dualizing complex $D$, a specific acyclic complex $K$, defined in
  \ref{K}, serves as a test complex for Gorensteinness in the
  following sense: The ring $R$ is Gorenstein if and only if
  $\tp{K}{D}$ is acyclic.  This is achieved by Theorem~\ref{thm1}. In
  general, $K$ is an acyclic complex of flat $R$-modules.
  Corollary~\ref{art} shows that if $R$ is an artinian local ring,
  then $K$ is a complex of projective modules, and
  \eqclbl{i}--\eqclbl{iii} above are equivalent with
  \begin{eqc}
  \setcounter{eqc}{3}
  \item The complex $K$ is totally acyclic.
  \end{eqc}
  
  Test complexes of injective modules can be obtained directly from
  $K$ (Corollary \ref{cor:K}) or through a potentially different
  construction explored in Section \ref{sec:section 3}. The authors of
  \cite{SInHKr06} have pointed out that the latter is of particular
  interest, as it yields a generator for $\Kac/\Ktac$, the Verdier
  quotient of acyclic complexes modulo totally acyclic complexes in
  the homotopy category of injective $R$-modules. This is proved in
  Theorem~\ref{thm3}.%
\end{ipg}

\section{Background}
\label{sec:1}

Throughout this paper $R$ is a commutative noetherian ring. The
notation $\Rmk$ means $R$ is local with maximal ideal $\m$ and residue
field $k$.

Complexes of $R$-modules ($R$\emph{-complexes} for short) are graded
homologically,
\begin{equation*}
  X = \dots \to X_{i+1} \xra{\partial_{i+1}^X} X_i
  \xra{\partial_{i}^X} X_{i-1} \to \cdots.
\end{equation*}
The \emph{suspension} of $X$ is denoted $\Susp{X}$; it is the complex
with $(\Susp{X})_i = X_{i-1}$ and dif\-ferential $\partial^{\Susp{X}} =
-\partial^X$. A complex $X$ is said to be \emph{bounded} if
$X_i=0$ for $|i| \gg 0$.

An isomorphism between $R$-complexes is denoted by a '$\is$'; we write
$X\is Y$ if there exists an isomorphism $X \xra{\is} Y$.  

A morphism between $R$-complexes is called a \emph{quasi-isomorphism},
and denoted $X \xre Y$ if the induced map in homology,
$\H{X}\to\H{Y}$, is an isomorphism.  Following \cite[sec.~1]{LLAHBF91}
we write $X\simeq Y$, if $X$ and $Y$ can be linked by a sequence of
quasi-isomorphisms with arrows in alternating directions.
Recall that a morphism $X\to Y$ is a quasi-isomorphism if and only if
its mapping cone, written $\Cone{(X\to Y)}$, is acyclic.

\begin{bfhpg}[Resolutions]
  \label{res}
  The following facts are established in
  \cite[sec.~1]{LLAHBF91}\footnote{ Where semi-projective/injective
    resolutions are called DG-projective/injective.} and \cite{dga}.
  
  Every $R$-complex $X$ has a semi-projective resolution. That is,
  there is a quasi-isomorphism $P \xre X$, where $P$ is a complex of
  projective $R$-modules such that $\Hom{P}{-}$ preserves
  quasi-isomorphisms. For such a complex, also the functor $\tp{-}{P}$
  preserves quasi-isomorphisms. In particular, for any $R$-complexes
  $Y\simeq Z$ we have $\Hom{P}{Y}\simeq \Hom{P}{Z}$ and
  $\tp{Y}{P}\simeq\tp{Z}{P}$.
  
  If there is an $l$ such that $\H[i]{X}=0$ for $i<l$, then $X$ has a
  semi-projective resolution $P$ with $P_i=0$ for $i<l$.  If, in
  addition, $\H[i]{X}$ is finitely generated for all $i$, then $P$ can
  be chosen with all modules $P_i$ finitely generated.
  
  Every $R$-complex $X$ has a semi-injective resolution. That is,
  there is a quasi-isomorphism $X \xre J$, where $J$ is a complex of
  injective $R$-modules such that $\Hom{-}{J}$ preserves
  quasi-isomorphisms. In particular, for such a complex $J$ and any
  $R$-complexes $Y\simeq Z$ we have $\Hom{Y}{J}\simeq \Hom{Z}{J}$.
\end{bfhpg}

\begin{lem}
  \label{tensor-lemma}
  Let $X$ and $Y$ be $R$-complexes such that either \mbox{$X_i=0$} for
  all $i \ll 0$ or $Y_i=0$ for all $i \gg 0$. If $\H{\tp{X_i}{Y}} = 0$
  for all $i \in \ZZ$, then $\H{\tp{X}{Y}} = 0$.
\end{lem}

\begin{prf}
  Let $E$ be a faithfully injective $R$-module. The complex
  $\tp{X}{Y}$ is acyclic if and only if $\Hom{\tp{X}{Y}}{E} \is
  \Hom{X}{\Hom{Y}{E}}$ is so. The claim is now immediate
  from~\cite[lem.~(2.4)]{CFH-}.
\end{prf}

\begin{bfhpg}[Dualizing complexes]
  \label{dc}
  Following \cite[V.\S2]{rad}, a \emph{dualizing complex} for $R$ is a
  bounded complex $D$ of injective $R$-modules such that $\H[i]{D}$ is
  finitely generated for all $i\in\ZZ$, and the homothety morphism
  \begin{equation*}
    \dmapdef{\chi^D}{R}{\Hom{D}{D}}
  \end{equation*}
  is a quasi-isomorphism. 
  
  Let $\Rmk$ be a local ring with a dualizing complex $D$. After
  suspensions we can assume $D$ is \emph{normalized},
  cf.~\cite[V.\S5]{rad}, in which case \cite[prop.~V.3.4]{rad} yields
  \begin{equation}
  \label{dcl}
    \H{\Hom{k}{D}} \is k.
  \end{equation}
  If $R$ is artinian, then $\E{k}$, the injective hull of the residue
  field, is a normalized dualizing complex for $R$.
\end{bfhpg}

\section{A test complex of flat modules}
\label{sec:section 2}

\begin{bfhpg}[A distinguished complex of flat modules]
  \label{K}
  Assume that $R$ has a dualizing complex $D$, and let $\pi\colon
  P\xre D$ be a semi-projective resolution. By \ref{res} we can assume
  that $P$ consists of finitely generated modules with $P_i=0$ for all
  $i\ll 0$. The functors $\Hom{P}{-}$ and $\Hom{-}{D}$ preserve
  quasi-isomorphisms, so the commutative
 diagram
  \begin{equation*}
    \xymatrixrowsep{2em}
    \xymatrixcolsep{5em}
    \xymatrix
      {
      \Hom{P}{P}
                \ar[r]^-{\Hom{P}{\pi}}_-{\eq}
      &\Hom{P}{D}\\
      R 
        \ar[r]^-{\chi^D}_-{\eq}
        \ar[u]^-{\chi^P}
      &\Hom{D}{D}
                \ar[u]_-{\Hom{\pi}{D}}^-{\eq}
    } 
  \end{equation*}
  shows that also the homothety map $\chi^P$ is a
  quasi-isomorphism. In particular,
  \begin{equation*}
   K = \Cone{(R \xra{\;\chi^P\;} \Hom{P}{P})}
  \end{equation*}
  is acyclic.  The modules in $\Hom{P}{P}$ are direct products of
  modules of the form $\Hom{P_i}{P_{i+n}}$, and each such module is
  flat.  Thus, $\chi^P$ is a quasi-isomorphism between complexes of
  flat $R$-modules, and the mapping cone $K$ is, therefore, an acyclic
  complex of flat $R$-modules.
\end{bfhpg}

We can now state the main result; the proof is given at the end of the
section.

\begin{thm}
  \label{thm1}
  Let $R$ be a commutative noetherian ring with a dualizing complex
  $D$, and let $K$ be the acyclic complex of flat modules defined
  in \ref{K}.  The ring $R$ is Gorenstein if and only if the complex
  $\tp{K}{D}$ is acyclic.
\end{thm}

\begin{rmk}
  While also $C=\Cone{\chi^D}$ is an acyclic complex of flat
  $R$-modules, it cannot detect Gorensteinness. Indeed, $C$ is
  bounded, so $\tp{C}{X}$ is acyclic for every $R$-complex $X$ by
  Lemma~\ref{tensor-lemma}. If $R$ is artinian, then $C$ is even split
  exact.
 \end{rmk}

\begin{rmk}
  \label{rmk:complete flat}
  In the theory of Gorenstein dimensions, there is a notion of a
  \emph{complete flat resolution}---due to Enochs, Jenda, and
  Torrecillas \cite{EJT-93}---namely an acyclic complex $F$ of flat
  modules such that $\tp{F}{I}$ is acyclic for every injective module
  $I$.
  
  If $R$ is Gorenstein, then every acyclic complex of flat $R$-modules
  is a complete flat resolution. Indeed, every injective $R$-module
  $I$ has finite flat dimension, and then it is straightforward to
  verify that the functor $\tp{-}{I}$ preserves acyclicity of
  complexes of flat modules.  On the other hand, let $K$ and $D$ be as
  in Theorem~\ref{thm1}. If $K$ is a complete flat resolution, then
  $\tp{K}{D}$ is acyclic by Lemma~\ref{tensor-lemma}.
  
  Thus, the following assertions are equivalent:
\begin{eqc}
 \item The ring $R$ is Gorenstein.
 \item The complex $K$ is a complete flat resolution.
 \item Every acyclic complex of flat modules is a complete flat
   resolution.
 \end{eqc}
\end{rmk}

\begin{ipg}
  The complex $K$ defined in \ref{K} appears to be a natural test
  object for Gorensteinness. However, it might in the context of
  \cite{SInHKr06} be of interest to exhibit a test complex of
  injective or of projective modules.
  
  To this end, we first note that the next corollary to
  Theorem~\ref{thm1} is immediate in view of Remark~\ref{rmk:complete
    flat} and \cite[prop.~(6.4.1)]{LWC}. See Section~\ref{sec:section
    3} for a further discussion of test complexes of injective
  modules.
\end{ipg}

\begin{cor}
  \label{cor:K}
  Let $R$ be a commutative noetherian ring with a dualizing complex.
  Let $K$ be the acyclic complex of flat modules defined in \ref{K},
  and let $E$ be a faithfully injective $R$-module. The complex
  $\Hom{K}{E}$ is an acyclic complex of injective modules, and $R$ is
  Gorenstein if and only if $\,\Hom{K}{E}$ is totally acyclic. \qed
\end{cor}

  For artinian local rings $(R,\m)$, Theorem~\ref{thm1} provides a
  test complex of projective modules. In particular, for $R$ with
  $\m^2=0$ the following recovers \cite[prop.~6.1(3)]{SInHKr06}.

\begin{cor}
  \label{art}
  Let $R$ be an artinian local ring. The complex $K$ defined in
  \ref{K} is an acyclic complex of projective $R$-modules, and $R$ is
  Gorenstein if and only if $K$ is totally acyclic.
\end{cor}

\begin{prf}
  When $R$ is artinian and local, every flat $R$-module is projective.
  Thus, $K$ is an acyclic complex of projective modules.
  
  The ``only if'' part is well-known. To prove ``if'', assume $K$ is
  totally acyclic and recall from \ref{dc} that the module $E=\E{k}$
  is dualizing for $R$. The first of the following isomorphisms is
  induced by $\chi^E$, and the second is Hom-tensor adjointness
  \begin{equation*}
    \Hom{K}{R} \is \Hom{K}{\Hom{E}{E}} \is \Hom{\tp{K}{E}}{E}.
  \end{equation*}
  The complex $\Hom{K}{R}$ is acyclic and $E$ is faithfully injective,
  so $\tp{K}{E}$ is acyclic and, therefore, $R$ is Gorenstein by
  Theorem~\ref{thm1}.
\end{prf}

  For the proof of Theorem~\ref{thm1} we need a technical lemma.

\begin{lem}
  \label{ev}
  Let $P$ be an $R$-complex of finitely generated projective modules,
  $X$ be any $R$-complex, and $B$ be a bounded $R$-complex of finitely
  generated modules. There is an isomorphism of $R$-complexes
  \begin{equation*}
    \dmapdef[\xri]{\omega}{\tp{\Hom{P}{X}}{B}}{\Hom{P}{\tp{X}{B}}}.
  \end{equation*}
\end{lem}

\begin{prf}
  It is straightforward to check that the assignment
  \begin{equation*}
    \omega(\phi\otimes b)(p) = (-1)^{|p||b|}\phi(p)\otimes b,
  \end{equation*}
  where $|\cdot |$ denotes the degree of an element, defines a morphism
  between the relevant complexes.  By assumption, there exist integers
  $l \le u$ such that $B_h =0$ when  $h<l$ or $h>u$.  For every
  $n\in\ZZ$ we have
  \begin{align*}
    \tpP{\Hom{P}{X}}{B}_n &= \bigoplus_{i=n-u}^{n-l}
    \tp{\Hom{P}{X}_i}{B_{n-i}}\\
    &= \bigoplus_{i=n-u}^{n-l}
    \tp{\big(\prod_{j\in\ZZ}\Hom{P_j}{X_{j+i}}\big)}{B_{n-i}}\\
    &\is \bigoplus_{i=n-u}^{n-l}
    \prod_{j\in\ZZ}\tpP{\Hom{P_j}{X_{j+i}}}{B_{n-i}}\\
    &\is \bigoplus_{i=n-u}^{n-l}
    \prod_{j\in\ZZ}\Hom{P_j}{\tp{X_{j+i}}{B_{n-i}}}\\
    &\is \prod_{j\in\ZZ}\Hom{P_j}{\bigoplus_{i=n-u}^{n-l}
      \tp{X_{j+i}}{B_{n-i}}}\\
    &= \prod_{j\in\ZZ}\Hom{P_j}{\tpP{X}{B}_{j+n}}\\
    &= \Hom{P}{\tp{X}{B}}_n.
  \end{align*}
  Since the modules $B_{n-i}$ are finitely generated, the functors
  $\tp{-}{B_{n-i}}$ commute with arbitrary products for every $i$; this
  explains the first isomorphism. The modules $P_j$ are finitely
  generated and projective, so
  for all $i$, $j$, and $n$ the homomorphism of modules
  \begin{equation*}
    \tp{\Hom{P_j}{X_{j+i}}}{B_{n-i}} \xra{\omega_{ijn}}
    \Hom{P_j}{\tp{X_{j+i}}{B_{n-i}}}
  \end{equation*}
  is invertible, and this accounts for the second isomorphism. Thus,
  $\omega$ is an isomorphism of graded modules, and the sign
  in the definition of $\omega$ ensures that it commutes with the differentials.
\end{prf}

\begin{prf}[Proof of Theorem \ref{thm1}]
  The ``only if'' part was settled in Remark \ref{rmk:complete flat}.
  
  For the ``if'' part, assume that the complex $\tp{K}{D}$ is acyclic;
  the isomorphism $\Cone{(\tp{\chi^P}{D})} \is \tp{K}{D}$ implies that
  \begin{equation*}
    \tag{1}
    \dmapdef{\tp{\chi^P}{D}}{D}{\tp{\Hom{P}{P}}{D}}
  \end{equation*}
  is a quasi-isomorphism. 
  
  Choose an $n$ such that $\H[i]{D}=0$ for all $i > n$, and let $B$ be
  the soft truncation of $P$ on the left at $n$:
  \begin{equation*}
    B = 0 \lora \Coker{\partial_{n+1}^{P}}
    \xra{\overline{\partial_{n}^{P}}} P_{n-1} \xra{\partial_{n-1}^P}
    P_{n-2} \lora \cdots.
  \end{equation*}
  There are quasi-isomorphisms $B \xle P \xre D$ and, hence, a
  quasi-isomorphism $\mapdef[\xre]{\beta}{B}{D}$; see \cite[1.1.I.(1)
  and 1.4.I]{LLAHBF91}.  Since the mapping cone of $\beta$ is a
  bounded acyclic complex, and $\Hom{P}{P}$ is a complex of flat
  modules, Lemma~\ref{tensor-lemma} applies to show that also
  $\tp{\Hom{P}{P}}{\Cone(\beta)}$ is acyclic.  Thus, the isomorphism
  $\Cone({\tp{\Hom{P}{P}}{\beta}})\is\tp{\Hom{P}{P}}{\Cone(\beta)}$
  implies that also
  \begin{equation*}
    \tag{2}
    \dmapdef{\tp{\Hom{P}{P}}{\beta}}{\tp{\Hom{P}{P}}{B}}{\tp{\Hom{P}{P}}{D}}
  \end{equation*}
  is a quasi-isomorphism. 
  
  By the choice of $P$, cf.~\ref{K}, the bounded complex $B$ consists
  of finitely generated modules, and Lemma~\ref{ev} yields an
  isomorphism
 \begin{equation*}
   \tag{3}
   \dmapdef[\xri]{\omega}{\tp{\Hom{P}{P}}{B}}{\Hom{P}{\tp{P}{B}}}.
 \end{equation*}

 Finally, let $\mapdef[\xra{\eq}]{\iota}{\tp{P}{B}}{J}$ be a semi-injective
 resolution; the quasi-isomorphism $\iota$ is preserved by
 $\Hom{P}{-}$, and the resulting quasi-isomorphism combined with (1),
 (2), and (3) yields
  \begin{equation*}
    \tag{4}
    D\simeq \Hom{P}{J}.
  \end{equation*}
  
  It suffices to prove that $R_\m$ is Gorenstein for every maximal
  ideal $\m$ of $R$. Let $\m$ be a maximal ideal; the complex $D_\m$
  is dualizing for $R_\m$, see \cite[cor.~V.2.3]{rad}. Set $k=R_\m/\m
  R_\m \is R/\m$. We may, after suspensions, assume $D_\m$ is
  normalized, so $\H{\Hom[R_\m]{k}{D_\m}} \is k$; see \pgref{dcl}.
  Moreover, there are isomorphisms $\Hom{k}{D} \is
  \tp{\Hom{k}{D}}{R_\m} \is \Hom[R_\m]{k}{D_\m}$, so we have
 \begin{equation*}
    \tag{5}
  k\is\H{\Hom{k}{D}}.
  \end{equation*}
  
  Let $\mapdef[\xre]{\upsilon}{Q}{k}$ be a semi-projective resolution of
  $k$ over $R$. As $\Hom{-}{D}$ preserves quasi-isomorphisms, we have
  \begin{equation*}
     \tag{6}
    \Hom{k}{D} \xre \Hom{Q}{D}.
  \end{equation*}
  Also $\Hom{Q}{-}$ preserves quasi-isomorphisms, and from (4) we get
  \begin{equation*}
     \tag{7}
     \Hom{Q}{D} \eq \Hom{Q}{\Hom{P}{J}}\is\Hom{\tp{Q}{P}}{J},
  \end{equation*}
  where the isomorphism is Hom-tensor adjointness. Finally,
  $\tp{\upsilon}{P}$ is a quasi-isomorphism, and hence so is
  \begin{equation*}
     \tag{8}
     \mapdef[\xre]{\Hom{\tp{\upsilon}{P}}{J}}{\Hom{\tp{k}{P}}{J}}
     {\Hom{\tp{Q}{P}}{J}}.
  \end{equation*}
  
  Combining $(5)$--$(8)$ and again using Hom-tensor adjointness, we
  obtain
  \begin{align*}
    k &\is\H{\Hom{\tp{k}{P}}{J}}\\
    &\is \H{\Hom{\tp[k]{\tpP{k}{P}}{k}}{J}}\\
    &\is\H{\Hom[k]{\tp{k}{P}}{\Hom{k}{J}}}\\
    &\is \Hom[k]{\H{\tp{k}{P}}}{\H{\Hom{k}{J}}}.
  \end{align*}
  Thus, $\Hom[k]{\H{\tp{k}{P}}}{\H{\Hom{k}{J}}}$ is a finitely
  generated $k$-vector space; in particular, $\H{\tp{k}{P}}$ must be
  finitely generated. Note that $\H[i]{\tp{k}{P}} \is
  \H[i]{\tp[R_\m]{k}{P_\m}}$ for all $i\in\ZZ$; it follows that
  $\H[i]{\tp[R_\m]{k}{P_\m}} =0$ for all $i\gg 0$. By
  \cite[prop.~5.5]{LLAHBF91} the dualizing $R_\m$-complex $D_\m$ then
  has finite flat dimension, and hence $R_\m$ is Gorenstein; see
  \cite[thm.~(17.23)]{hha} or \cite[thm.~(8.1)]{LWC01a}.
\end{prf}
 
\section{A test complex of injective modules}
\label{sec:section 3}

The next construction is another source for test complexes.

\begin{bfhpg}[A distinguished complex of injective modules]
  \label{M}
  Assume $R$ has a dualizing complex $D$. As in \ref{K}, let
  $\mapdef[\xre]{\pi}{P}{D}$ be a semi-projective resolution of $D$
  consisting of finitely generated modules with $P_i=0$ for all $i\ll
  0$. The assignment $\varphi \otimes p \mapsto \varphi(p)$ defines a
  morphism of complexes, $\varepsilon$, such that the following
  diagram is commutative:
  \begin{equation*}
    \xymatrix@C=4em{\tp{\Hom{P}{D}}{P} \ar[r]^-{\varepsilon} & D\\
    \tp{\Hom{P}{P}}{P} \ar[u]_-{\eq}^-{\tp{\Hom{P}{\pi}}{P}} & 
    \ar[l]_-{\eq}^-{\tp{\chi^P}{P}} \ar[u]^{\eq}_-{\pi} \tp{R}{P}.}
  \end{equation*}
  Thus, $\varepsilon$ is a quasi-isomorphism between complexes of
  injective $R$-modules, and the mapping cone
  \begin{equation*}
    \label{eq:M}
    M = \Cone{(\tp{\Hom{P}{D}}{P} \xra{\;\varepsilon\;} D)}
  \end{equation*}
   an acyclic complex of injective $R$-modules.
\end{bfhpg}

  An argument similar to the proof of Theorem \ref{thm1} yields the
  next result, which is also a corollary of Theorem~\ref{thm3}.

\begin{thm}
  \label{thm2}
  Let $R$ be a commutative noetherian ring with a dualizing complex,
  and let $M$ be the acyclic complex of injective modules defined in
  \ref{M}.  The ring $R$ is Gorenstein if and only if $M$ is totally
  acyclic. \qed
\end{thm}

\begin{rmk}
  If $\Rmk$ is an artinian local ring, then there is an isomorphism
  \begin{equation*}
    K \is \Susp{\Hom{M}{\E{k}}}
  \end{equation*}
  where $K$ and $M$ are the complexes from \ref{K} and \ref{M}, and
  $\E{k}$ is the injective hull of $k$. Indeed, with $E=\E{k}$ there
  is a commutative diagram
  \begin{equation*}
    \xymatrixrowsep{2em}
    \xymatrixcolsep{4em}
    \xymatrix{
      \Hom{E}{E}\ar[r]^-{\Hom{\varepsilon}{E}}_-{\eq}
      \ar@{<-}[dd]_-{\cong}^-{\chi^E}
      & \Hom{\tp{\Hom{P}{E}}{P}}{E} \\
      & \ar[u]^{\cong} \Hom{P}{\Hom{\Hom{P}{E}}{E}}\\
      R \ar[r]_-{\eq}^-{\chi^P}&\ar[u]^{\cong} \Hom{P}{P}.
    } 
  \end{equation*}
  The vertical maps on the right are the natural isomorphisms, and
  because $E$ is a module, also the homothety map $\chi^E$ is a
  genuine isomorphism.  The diagram induces the desired isomorphism
  between the complexes $K = \Cone{(\chi^P)}$ and
  $\Cone{(\Hom{\varepsilon}{E})} \is
  \Susp{\Hom{\Cone{(\varepsilon)}}{E}} = \Susp{\Hom{M}{E}}$.
  
  When $R$ is not artinian, we do not know if the complexes $K$ and
  $M$ are related.
\end{rmk}

  Using \cite[prop.~(5.1)]{CFH-} it is not hard to prove the
  next parallel to Corollary~\ref{cor:K}.

\begin{cor}
  \label{cor:M}
  Let $R$ be a commutative noetherian ring with a dualizing complex.
  Let $M$ be the acyclic complex of injective modules defined in
  \ref{M}, and let $E$ be a faithfully injective $R$-module. The complex
  $\Hom{M}{E}$ is an acyclic complex of flat modules, and $R$ is
  Gorenstein if and only if $\,\tp{\Hom{M}{E}}{I}$ is acyclic for every
  injective module $I$. \qed
\end{cor}

In conversations, the authors of \cite{SInHKr06} have informed us of
Theorem ~\ref{thm3} below; note that it contains Theorem~\ref{thm2}.
For notation and terminology we refer to \cite{SInHKr06}.
  
\begin{thm} 
  \label{thm3}
  Let $R$ be a commutative noetherian ring with a dualizing complex.
  The acyclic complex $M$ of injective modules defined in \ref{M}
  generates the quotient category $\Kac/\Ktac$.
\end{thm}

\begin{proof}
  By \cite[1.7, 5.4, and 5.9(3)]{SInHKr06} the quotient category
  $\Kac/\Ktac$ is generated by the image of the dualizing complex $D$
  under the equivalence\linebreak \mbox{$\Df \xra{\;\sim\;} \Kcinj$},
  cf.~\cite[2.3(2)]{SInHKr06}.
  
  Let $P \xre D$ be a semi-projective resolution. The functor
  $\Hom{P}{-}$ preserves quasi-isomorphisms, so the composite
  \begin{equation*}
    R \xre \Hom{P}{P} \xre \Hom{P}{D}
  \end{equation*}
  provides an injective resolution $R \xre \iR = \Hom{P}{D}$. Since
  $\iR$ is a compact object in $\Kinj$, the inclusion of the
  localizing subcategory $\Loc(\iR) \subseteq \Kinj$ admits a right
  adjoint $\mapdef{\rho}{\Kinj}{\Loc(\iR)}$; see~\cite[1.5.1]{SInHKr06}.
  By \cite[2.3(2)]{SInHKr06} the image of $D$ in $\Kinj$ is
  \begin{equation*}
    \Cone(\rho(D)\xra{\xi} D),
  \end{equation*}
  where $\xi$ is the natural map.
  
  It remains to show that $M\cong \Cone(\rho(D)\xra{\xi} D)$. It
  suffices to establish a commutative diagram,
  \begin{equation*}
    \xymatrix{ 
      \rho(D)\ar[r]^-{\xi}\ar[d]^-{\cong}&D\\
      \tp{\Hom{P}{D}}{P}\ar[ur]^-{\varepsilon}&{}\\
    }
  \end{equation*}
  The complex $\tp{\Hom{P}{D}}{P} = \tp{\iR}{P}$ is in $\Loc(\iR)$, and
  since $\varepsilon$ is a quasi-isomorphism,
  $\Hom[\Kinj]{\iR}{\varepsilon}$ is an isomorphism,
  cf.~\cite[2.2]{SInHKr06}. The existence of the desired isomorphism
  $\rho(D) \is \tp{\Hom{P}{D}}{P}$ now follows from
  \cite[1.4]{SInHKr06}.
\end{proof}


\section*{Acknowledgments}
It is a pleasure to thank Srikanth Iyengar and Henning Krause for
discussions regarding this work and, in particular, for letting us
include Theorem~\ref{thm3}. Thanks are also due to Paul Roberts
for suggestions that have improved the exposition.



\providecommand{\bysame}{\leavevmode\hbox to3em{\hrulefill}\thinspace}
\renewcommand{\MR}{\relax\ifhmode\unskip\space\fi MR }
\providecommand{\href}[2]{#2}

\end{document}